\documentclass[12pt]{amsart}

\usepackage{amsmath}
\usepackage{amssymb}
\newtheorem{theorem}{Theorem}[section]

\newtheorem{corollary}[theorem]{Corollary}

\theoremstyle{definition}

\begin{document}

\title{On the  Transcendence of  certain   real numbers}
\author{Veekesh Kumar}
\address{Harish-Chandra Research Institute, HBNI\\
Chhatnag Road, Jhunsi\\
Allahabad - 211019\\
INDIA}
\email[Veekesh Kumar]{veekeshkumar@hri.res.in}

\author[B. Mance]{Bill Mance}
\address[B. Mance]{Uniwersytet im. Adama Mickiewicza w Poznaniu, Collegium Mathematicum, ul. Umultowska 87, 61-614 Pozna\'{n}, Poland}

\email{Bill.A.Mance@gmail.com}

\subjclass[2010] {Primary 11J68; Secondary 11K16}
\keywords{Rational approximation, Subspace theorem.}
\bigskip

\begin{abstract} 
In this article we discuss the transcendence of certain infinite sums and products by using the Subspace theorem. In particular we improve the result of Han\v{c}l and Rucki \cite{hancl3}.
\end{abstract}
\maketitle

\section{Introduction}
 There are several methods  to prove the transcedence of an infinite series.  By using Mahler's method \cite{nish} one can prove transcendence of certain infinite sums.  In 2001, Adhikari {\it et al}. \cite{shorey}  studied the transcendence of certain infinite series by using a theorem of Baker regarding linear forms in logarithms of algebraic numbers. In the same year,  Han\v{c}l \cite{hancl1} and Nyblom \cite{nyb} (see also \cite{nyb1}) studied the transcendence of infinite series using Roth's Theorem. In $2004$,  Adamczewski, Bugeaud, and Luca \cite{adam2} (see also \cite{adam1}) proved  a criterion for a $b$-ary expansion of real number to be transcendental, using the Subspace theorem. 
\bigskip

We state the following     conditions on the  sequences $(c_n)_n$, $(c'_n)_n$ and $(b_n)_n$  of positive integers.
\bigskip

\noindent{\bf Condition 1. }   Let $\delta> 0$ be a real number.  For any given integer $m\geq 2$,  let $(c_{i,n})_n$, $i=1,2,\ldots,m$ be a  sequence of non-zero integers and $(b_n)_n$ be the sequence of positive integers such that
\begin{equation*}\label{eq1.1}
\tag{1.1}
\limsup_{n\rightarrow\infty}\frac{b_{n+1}}{(b_1 b_2\ldots b_n)^{1+\delta}}\frac{1}{c_{i,n+1}}=\infty
\end{equation*}
and 
\begin{equation*}\label{eq1.2}
\tag{1.2}
\liminf_{n\rightarrow\infty}\frac{b_{n+1}}{b_n}\frac{c_{i,n}}{c_{i,n+1}}>1 \quad \mbox{~~for all ~~}\quad i\in\{1,2,\ldots,m\}.
\end{equation*}

\noindent{\bf Condition 2.}  Let $\delta$ and $\epsilon$ be positive real numbers. Let $\delta> 0$ be a real number.  For any given integer $m\geq 2$,  let $(c_{i,n})_n$, $i=1,2,\ldots,m$, be a  sequence of non-zero integers and let $(b_n)_n$ be a sequence of positive integers such that
\begin{equation*}\label{eq1.3}
\tag{1.3}
\limsup_{n\rightarrow\infty}\frac{b_{n+1}}{(b_1 b_2\ldots b_n)^{1+\delta+\frac{1}{\epsilon}}}\frac{1}{c_{i,n+1}}=\infty
\end{equation*}
and for all sufficiently large values of $n$
\begin{equation*}\label{eq1.4}
\tag{1.4}
\sqrt[1+\epsilon]{\frac{b_{n+1}}{c_{i,n+1}}}\geq \sqrt[1+\epsilon]{\frac{b_n}{c_{i,n}}}+1 \quad \mbox{~~for all ~~}\quad i\in\{1,2,\ldots,m\}.
\end{equation*} 
  
 In $2005$,  J. Han\v{c}l and Rucki \cite{hancl3} gave sufficient conditions under which an infinite sum is transcendental.  We mention one of their results here. 
\bigskip
 
\noindent{\bf Theorem A.} {\it Let $\delta>0$ be a real number.  Let $(c_n)_n$ and $(b_n)_n$ be  sequences of positive integers such that
\begin{equation*}
\limsup_{n\rightarrow\infty}\frac{b_{n+1}}{(b_1 b_2\ldots b_n)^{2+\delta}}\frac{1}{c_{n+1}}=\infty
\end{equation*}
and 
\begin{equation*}
\liminf_{n\rightarrow\infty}\frac{b_{n+1}}{b_n}\frac{c_n}{c_{n+1}}>1.
\end{equation*}
Then the real number 
$$
\alpha=\sum_{n=1}^\infty\frac{c_n}{b_n}
$$
is transcendental.}
\bigskip

In this article we give an improvement of the results in \cite{hancl3}, in particular an improvement of Theorem A,  and we also  study the transcendence of infinite products.  We now state our main result.
\begin{theorem}\label{maintheorem}
For any given integer $m\geq 2$, let  $\alpha_1, \alpha_2,\ldots,\alpha _m $ be real numbers.  Let  $\delta> 0$ be a real number such that  $\delta>\frac{1}{m}$.    Suppose there exist infinitely many $(m+1)$-tuples $(p_{n1}, p_{n2}, \ldots p_{nm}, q_n)$  of integers satisfying $q_n\neq 0$  and 
\begin{equation*}\label{eq1.5}
\tag{1.5}
\left|\alpha_i -\frac{p_{in}}{q_n}\right|<\frac{1}{q_n^{1+\delta}}, \quad \mbox{~~~for~~} 1\leq i\leq r.
\end{equation*}
Then either  the real numbers $1, \alpha_1, \alpha_2, \ldots,\alpha _m $  are $\mathbb{Q}$- linearly dependentat or  at least one of $\alpha_i$'s is transcendental. 
\end{theorem} 
As an application of Theorem \ref{maintheorem}, we have the following results about transcendence.
\begin{theorem}\label{maintheorem1}
Let  $Q=(b_n)_{n\geq 1}$ be  a sequence of positive integers   with $b_n\geq 2$ for all integers $n\geq 1$ and let $\delta>\frac{1}{3}$ be a real number.   Suppose that there exists  an infinite subset $T$   of natural numbers  $N$   such that
\begin{equation*}\label{eq1.6}
\tag{1.6}
 \sigma(N+1)(b_1b_2 \ldots b_N)^\delta \leq b_{N+1};
\end{equation*} 
\begin{equation*}\label{eq1.7}
\tag{1.7}
\phi(N+1)(b_1b_2 \ldots b_N)^\delta \leq b_{N+1};
\end{equation*}
and 
\begin{equation*}\label{eq1.8}
\tag{1.8}
d(N+1)(b_1b_2 \ldots b_N)^\delta \leq b_{N+1}.
\end{equation*}
Then at least one of the real numbers
$$
\beta_1=\sum_{n=1}^\infty\frac{\sigma(n)}{b_1 b_2\ldots b_n} ,\quad  \beta_2=\sum_{n=1}^\infty\frac{\phi(n)}{b_1 b_2\ldots b_n},\quad
\beta_3=\sum_{n=1}^\infty\frac{d_n}{b_1 b_2\ldots b_n}
$$
is transcendental.
\end{theorem}
\begin{theorem}\label{maintheorem2}
For any given integer $m\geq 2$,  let $\delta>\frac{1}{m}$ be a real number. Let $(c_{i,n})_n$, $i=1,2,\ldots,m$  and $(b_n)_n$ be  sequences of positive integers satisfying {\it Condition 1}.  
Then either at least one of the real numbers
$$
\beta_1=\sum_{n=1}^\infty\frac{c_{1,n}}{b_n},\quad \beta_2=\sum_{n=1}^\infty\frac{c_{2,n}}{b_n},  \cdots, \beta_m=\sum_{n=1}^\infty\frac{c_{m,n}}{b_n}
$$
is transcendental or $1, \beta_1, \beta_2,\ldots, \beta_m$ are $\mathbb{Q}$- linearly dependent.
\end{theorem}
The following corollary  of Theorem \ref{maintheorem2} tells us that Theorem \ref{maintheorem2} is an improvement of Theorem A. 
\begin{corollary}\label{cor}   Let $\delta>\frac{1}{2}$ be real number and let $(b_n)_n$ be a sequence of positive integers such that $b_1=2 $ and 
$$
 b_{n+1}=(b_1 b_2 \cdots b_n)^2, \quad\mbox{for all integers } n\geq 2.
$$
Then at least one of the real numbers
$$
\sum_{n=1}^\infty\frac{1}{b_n} \quad \mbox{and} \quad \sum_{n=1}^\infty\frac{d(n)}{b_n}
$$
is transcendental, where $d(n)=\sum_{d|n}1$.
\end{corollary}   
First we  note that $(b_1 b_2 \cdots b_n)^{1+\delta}\leq b_{n+1}$ for any $\frac{1}{2}<\delta<1$.   But $(b_1 b_2 \cdots b_n)^{2+\delta}> b_{n+1}$ for any choice of $\delta>0$. Therefore by Theorem A  we can not conclude the transcendence of  any of the numbers  
$$
 \sum_{n=1}^\infty\frac{c_{1,n}}{b_n} \quad \mbox{and} \quad \sum_{n=1}^\infty\frac{c_{2,n}}{b_n}.
$$
Theorem \ref{maintheorem2} can be improved by saying that at least one of $\beta_i$ is transcendental under some more assumptions on the sequences $(c_{i,n})_n$ and $(b_n)_n$.  More precisely, we have the following theorem.
\begin{theorem}\label{maintheorem3}
 For any given integer $m\geq 2$, let $\delta>\frac{1}{m}$ be a real number.  Let $(c_{i,n})_n$, $i=1,2,\ldots,m$  and $(b_n)_n$ be  sequences of positive integers satisfying {\it Condition 1}. 
If we assume that 
$$
1\leq \liminf_{n\to\infty}b_n^{\frac{1}{(m+1)^n}}<\limsup_{n\to\infty}b_n^{\frac{1}{(m+1)^n}}<\infty
$$
and
$$
\lim_{n\to\infty}\frac{c_{i,n}}{c_{j,n}}=0,\quad \mbox{~~for  all~~~}\quad i,j\in\{1,2,\ldots,m\}, i>j, 
$$
then at least one of the real numbers
$$
\beta_1=\sum_{n=1}^\infty\frac{c_{1,n}}{b_n},\quad \beta_2=\sum_{n=1}^\infty\frac{c_{2,n}}{b_n},  \cdots, \beta_m=\sum_{n=1}^\infty\frac{c_{m,n}}{b_n}
$$
is transcendental.
\end{theorem}
\begin{theorem}\label{maintheorem4}
For any given integer $m\geq 2$, let $\delta$ and $\epsilon$ be positive real numbers such that $\frac{\delta \epsilon}{1+\epsilon}>\frac{1}{m}$.  Let $(c_{i,n})_n$, $i=1,2,\ldots,m$  and $(b_n)_n$ be  sequences of positive integers satisfying {\it Condition 2}. Then  at least one of the real numbers
$$
\beta_1=\sum_{n=1}^\infty\frac{c_{1,n}}{b_n},\quad \beta_2=\sum_{n=1}^\infty\frac{c_{2,n}}{b_n},  \cdots, \beta_m=\sum_{n=1}^\infty\frac{c_{m,n}}{b_n}
$$
is transcendental or $1, \beta_1, \beta_2,\ldots, \beta_m$ are $\mathbb{Q}$- linearly dependent. 
\end{theorem}
Theorem \ref{maintheorem4} can be improved by saying at least one of $\beta_i$ is transcendental under some more assumptions on the sequences $(c_{i,n})_n$ and $(b_n)_n$.  More precisely, we have the following theorem.
\begin{theorem}\label{maintheorem5}
For any given integer $m\geq 2$, let $\delta$ and $\epsilon$ be positive real numbers such that $\frac{\delta \epsilon}{1+\epsilon}>\frac{1}{m}$.  Let $(c_{i,n})_n$, $i=1,2,\ldots,m$  and $(b_n)_n$ be  sequences of positive integers satisfying {\it Condition 2}.     
If we assume that 
$$
1\leq \liminf_{n\to\infty}b_n^{\frac{1}{(m+1)^n}}<\limsup_{n\to\infty}b_n^{\frac{1}{(m+1)^n}}<\infty
$$
and
$$
\lim_{n\to\infty}\frac{c_{i,n}}{c_{j,n}}=0,\quad \mbox{~~for  all~~~}\quad i,j\in\{1,2,\ldots,m\}, i>j, 
$$
then at least one of the real numbers
$$
\beta_1=\sum_{n=1}^\infty\frac{c_{1,n}}{b_n},\quad \beta_2=\sum_{n=1}^\infty\frac{c_{2,n}}{b_n},  \cdots,\beta_m=\sum_{n=1}^\infty\frac{c_{m,n}}{b_n}
$$
is transcendental.
\end{theorem}
\begin{theorem}\label{maintheorem6}
For any given integer $m\geq 2$, let $\delta>\frac{1}{m}$ be a real number. Let $(c_{i,n})_n$,  $i=1,2,\ldots,m$ and $(b_n)_n$  be sequences of positive integers satisfying all the hypotheses of Theorem \ref{maintheorem2}. Suppose that $c_{i,n}\leq b_n$ for all $n\geq 1$ and  $i=1,2,\ldots,m$. Then at least one of the real numbers 
$$
\beta_1=\prod_{n=1}^\infty\left(1+\frac{c_{1,n}}{b_n}\right), \quad \beta_2=\prod_{n=1}^\infty\left(1+\frac{c_{2,n}}{b_n}\right),\cdots,\beta_m=\prod_{n=1}^\infty\left(1+\frac{c_{m,n}}{b_n}\right) 
$$
is transcendental. 
\end{theorem}

\section{preliminaries}
The Subspace theorem is the main tool to prove our results.  
In order to state the Subspace theorem, we need to introduce $p$-adic absolute values on a finite extension $K$ over ${\mathbb Q}$.  First, we shall define the $p$-adic absolute value on $\mathbb{Q}$ and then we shall extend this absolute value for a finite extension $K$ over ${\mathbb{Q}}$. 
\bigskip

Let $p$ be a prime number in $\mathbb{Z}$. Let $x/y$ be any rational number  where $x\in {\mathbb{Z}}\backslash\{0\}$, $y\geq 1$ integer and  $(x,y) = 1$. We define 
$$\mbox{ord}_p(x/y) = \left\{\begin{array}{ll}
 n &  \mbox{ if } p^n\Vert x\\
 -n &\mbox{ if } p^n\Vert y.
 \end{array}\right.
 $$ 
Then, the $p$-adic absolute value on $\mathbb{Q}$, denoted by $| \cdot |_p$ and defined  as 
$$
\left|\frac{x}{y}\right|_p = \left(\frac{1}{p}\right)^{\mbox{ord}_p(x/y)} \mbox{ and }  |0|_p = 0.$$
In this set up, the usual absolute value  $|\cdot|$ on $\mathbb{Q}$ is denoted by $| \cdot |_{\infty}$.  

\smallskip

Now, let $K/{\mathbb{Q}}$ be a number field and ${\mathcal{O}}_K$ be its ring of integers.  
Then, for any prime number $p\in \mathbb{Z}$,  the ideal $p\mathcal{O}_K$ in $\mathcal{O}_K$ can be factored into product of prime ideals as 
$$
p\mathcal{O}_K =  {\mathfrak{p}}_1^{e_1} \ldots \mathfrak{p}_g^{e_g}$$
with $e_i\geq 1$ integers and ${\mathfrak p}_i$ are prime ideals in $\mathcal{O}_K$. 
Hence $p{\mathcal O}_K \subset {\mathfrak p}_i$ for all $i = 1, 2, \ldots, g$. In this situation, we say ${\mathfrak p}_i | p$ (${\mathfrak p}_i$ divides $p$) for all $i = 1, 2, \ldots, g$. 

\bigskip
 
Since $K$ is the quotient field of ${\mathcal{O}}_K$,   any $\alpha \in K$ can be written as  
$\alpha = x/y$ where $x, y\in {\mathcal{O}}_K$ with gcd$(x{\mathcal{O}}_K, y{\mathcal{O}}_K) = {\mathcal{O}}_K$. Therefore, for any $\alpha\in K$ and  for a given  prime ideal $\mathfrak{p}$ in ${\mathcal{O}}_K$, we can define 
$$
\mbox{ord}_{\mathfrak{p}}(\alpha) = \left\{\begin{array}{ll}
n &\mbox{ if } {\mathfrak{p}}^n\Vert x{\mathcal{O}}_K\\
-n &  \mbox{ if } {\mathfrak{p}}^n\Vert y{\mathcal{O}}_K\\
\end{array} \right.
$$
Also, for any non-zero prime ideal $\mathfrak{p}$ in ${\mathcal{O}}_K$, the norm of ${\mathfrak{p}}$ denoted by $N{\mathfrak{p}}$ and defined by $N{\mathfrak{p}} = |{\mathcal{O}}_K/{\mathfrak{p}}|$, cardinality of the quotient ring (which is known to be finite).  Now, we can extend the $p$-adic absolute value for any $\alpha \in K\backslash\{0\}$ as 
$$
|\alpha|_p = \left|\frac{x}{y}\right|_p = \prod_{{\mathfrak{p}}|p} \left(\frac{1}{N{\mathfrak{p}}}\right)^{\mbox{ord}_{\mathfrak p}(\alpha)}.$$
If $p = \infty$, then we define 
$$
|\alpha|_\infty = |N_{K/{\mathbb{Q}}}(\alpha)|,$$
where $N_{K/{\mathbb{Q}}}(\alpha)$ is the norm of $\alpha$ (which is nothing but the product of all the Galois conjugates of $\alpha$) in $K/\mathbb{Q}$.  With these definitions, one can check the product formula
$$
|\alpha|_{\infty}\prod_{p}|\alpha|_p = 1
$$
holds for all $\alpha \in K\backslash\{0\}$.
\bigskip

For any vector $(y_1, \ldots, y_n) \in \mathbb{Z}^n\backslash (0,0,\ldots, 0)$, we define the {\it height}
$$
H((y_1, \ldots, y_n)) = \prod_{p} \max_{1\leq i\leq n} \{1,|y_1|_p,\ldots, |y_n|_p\}=
\max_{1\leq i\leq n} \{|y_1|_\infty,\ldots, |y_n|_\infty\}.
$$
We are now ready to state the $p$-adic version of the Subspace theorem 

\begin{theorem} \label{subspace} (Subspace theorem)  
Let $S_f$ be a finite subset of prime numbers in $\mathbb{Z}$ and let $S = S_f \cup \{\infty\}$. Let $n> 1$ be an integer. For every prime $p \in S$, let the given linear forms $L_{1, p}, \ldots, L_{n,p}$ in $n$-variables, whose coefficients are algebraic numbers,  be  linearly independent.   Let  $\epsilon > 0$ be given. Then    the set 
$$
T=\left\{(y_1, \ldots, y_n)\in \mathbb{Z}^{n}: \prod_{p\in S} \prod_{i=1}^n|L_{i,p}{(y_1, \ldots, y_n)}|_p \leq H((y_1, \ldots, y_n)) ^{-\epsilon}\right\}
$$
is contained in a finite union of proper subspaces of ${\mathbb{Q}}^n$. 
\end{theorem}
We need the following result of Erd\H {o}s and Strauss \cite{erdos}.
\begin{theorem}\label{erdos}
Let  $Q=(b_n)_{n\geq 1}$ be  a sequence of positive integers   with $b_n\geq 2$ for all integers $n\geq 1$ and let $\delta>\frac{1}{3}$ be any positive real number.   Suppose there exists  an infinite subset $T$  consists of natural numbers  $N$   such that
\begin{equation*}
  \sigma(N+1)(b_1b_2 \ldots b_N)^\delta \leq b_{N+1};
\end{equation*} 
\begin{equation*}
\phi(N+1)(b_1b_2 \ldots b_N)^\delta \leq b_{N+1}.
\end{equation*}
and 
\begin{equation*}
d(N+1)(b_1b_2 \ldots b_N)^\delta \leq b_{N+1}.
\end{equation*}
Then  the real numbers
$$
1,\quad \sum_{n=1}^\infty\frac{\sigma(n)}{b_1 b_2\cdots b_n} ,\quad  \sum_{n=1}^\infty \frac{\phi(n)}{b_1 b_2 \cdots b_n}, \mbox{~~and~~}  \alpha''=\sum_{n=1}^\infty\frac{d_n}{b_1 b_2\ldots b_n}
$$
are $\mathbb{Q}$-linearly independent,  where $\phi(n)$ denotes the Euler totient function,  $\displaystyle\sigma(n) = \sum_{d|n} d$ and $(d_n)_n$ is any sequence of integers satisfying $|d_n|<n^{\frac{1}{2}-\delta}$ for all large $n$ and $d_n \neq 0$ for infinitely many values $n$.
\end{theorem}  
We need the following result of Han\v{c}l \cite{hancl4}
\begin{theorem}\label{maintheorem7}
For given integer $m\geq 2$, let $(c_{i,n})_n$, $i=1,2,\ldots,m$ be a  sequence of positive integers integers and let $(b_n)_n$ be a sequence of positive integers such that 
$$
1\leq \liminf_{n\to\infty}b_n^{\frac{1}{(m+1)^n}}<\limsup_{n\to\infty}b_n^{\frac{1}{(m+1)^n}}<\infty
$$
and
$$
\lim_{n\to\infty}\frac{c_{i,n}}{c_{j,n}}=0,\quad \mbox{~~for  all~~~}\quad i,j\in\{1,2,\ldots,m\}, i>j.
$$
Suppose that for all sufficiently large values of $n$,
$$
c_{i,n}<2^{(\log b_n)^\alpha},\quad i=1,2,\ldots,m,
$$
$$
b_n\geq n^{1+\epsilon}.
$$
Then the real numbers  
$$
1, \quad \sum_{n=1}^\infty\frac{c_{1,n}}{b_n},\cdots, \sum_{n=1}^\infty\frac{c_{m,n}}{b_n}
$$
are $\mathbb{Q}$-linearly independent.
\end{theorem}

We need the following result due to Han\v{c}l, Kolouch and Novotn\' y \cite{hancl5}.
\begin{theorem}\label{maintheorem8}
For any given integer $m\geq 2$, let $(c_{i,n})_n$, $i=1,2,\ldots,m$ be a  sequence of positive integers integers and let $(b_n)_n$ be a sequence of positive integers such that 
$$
1\leq \liminf_{n\to\infty}b_n^{\frac{1}{(m+1)^n}}<\limsup_{n\to\infty}b_n^{\frac{1}{(m+1)^n}}<\infty
$$
and
$$
\lim_{n\to\infty}\frac{c_{i,n}}{c_{j,n}}=0,\quad \mbox{~~for  all~~~}\quad j,i\in\{1,2,\ldots,m\}, i>j.
$$
Suppose that for all sufficiently large values of $n$,
$$
c_{i,n}<b^{\frac{1}{\log^{1+\epsilon}\log b_n}}_n,\quad i=1,2,\ldots,m,
$$
$$
b_n\geq n^{1+\epsilon}.
$$
Then the real numbers  
$$
1, \quad \prod_{n=1}^\infty\left(1+\frac{c_{1,n}}{b_n}\right),\cdots, \prod_{n=1}^\infty\left(1+\frac{c_{m,n}}{b_n}\right)
$$
are $\mathbb{Q}$-linearly independent.
\end{theorem}
\section{Proof of Theorem \ref{maintheorem}}
Without loss of generality we can assume that the absolute value of the real numbers  $\alpha_1, \alpha_2, \ldots, \alpha_m$ is less than $1$.  By \eqref{eq1.5}, we have 
\begin{align*}\label{eq3.1}
|q\alpha_1 - p_{n1}|&<\frac{1}{q_n^{\delta}};\\
|q\alpha_2 - p_{n2}|&<\frac{1}{q_n^{\delta}}\\
\vdots\\
|q\alpha_m - p_{nm}|&<\frac{1}{q_n^{\delta}}.\tag{3.1}
\end{align*}

Suppose the real numbers  $\alpha_1, \alpha_2, \ldots\alpha_m$ are  algebraic. In order to finish the proof we prove that  $1, \alpha_1, \alpha_2, \ldots,\alpha _m $  are $\mathbb{Q}$- linearly independent.
\bigskip

In order to prove  $1, \alpha_1, \alpha_2, \ldots,\alpha _m $  are $\mathbb{Q}$- linearly independent. We shall apply  Theorem \ref{subspace}.  Let  $S=\{\infty\}$ and  consider linear forms with algebraic coefficients
\begin{align*}\label{eq3.2}
L_{1, \infty}(X_1, X_2,\ldots, X_m) &= X_1,\\  L_{2, \infty}(X_1, X_2,\ldots, X_m) &= \alpha_1 X_1 - X_2,\\ 
\vdots\\
L_{m, \infty}(X_1, X_2,\ldots ,X_m) &= \alpha_r X_1 - X_m. \tag{3.2}
\end{align*}
Clearly, the  above linear forms  are linearly independent.
\bigskip

To apply   Theorem \ref{subspace}, we need to  compute the quantity
\begin{equation*}
\prod_{p\in S}  \prod_{i=1}^{m+1}|L_{i,p}{(p_{n1}, p_{n2}, \ldots,p_{nm}, q_n)}|_p=\prod_{i=1}^{m+1}|L_{i,\infty}{(p_{n1}, p_{n2},\ldots,p_{nm},q_n)}|_\infty.
\end{equation*}
First, by \eqref{eq3.2},  we note that
\begin{equation*}\label{eq3.3}
\tag{3.3}
\prod_{i=1}^{m+1}|L_{i,\infty}{((p_{n1}, p_{n2},\ldots,p_{nm},q_n))}|_\infty= |q_n|_\infty \prod_{i=1}^m|q_n \alpha_i-p_i|_\infty.
\end{equation*}
 Thus, from \eqref{eq3.1},  we conclude that
 $$
 \prod_{p\in S}  \prod_{i=1}^{m+1}|L_{i,p}{(p_{n1}, p_{n2}, \ldots,p_{nm}, q_n)}|_p<(H((p_{n1}, p_{n2}, \ldots,p_{nm},  q_n)))^{-\delta'}
 $$
holds for infinitely many  integers $n$  and for some $\delta'>0$.   Hence, by  Theorem \ref{subspace}, all these   non-zero integer lattice points  $X^{(n)}=(p_{n1}, p_{n2}, \ldots,p_{nm}, q_n)$  lie only in  finitely many proper subspaces of $\mathbb{Q}^{m+1}$.  Therefore,  there exists a proper subspace of $\mathbb{Q}^{m+1}$ containing these  integer lattice points. That is,  there exists a non-zero tuple $(z_1,z_2,\ldots, z_{m+1})\in \mathbb{Z}^{m+1}$ such that 
$$
z_1 q_n+z_2 p_{n1}+z_3 p_{n2}+\cdots+z_{m+1} p_{nm}=0
$$
holds for infinitely many  values of $n$.  This implies that 
\begin{align*}
\lim_{n\to\infty}\left(z_1+z_2\frac{p_{n1}}{q_n}+z_3\frac{p_{n2}}{q_n}+\cdots+z_{m+1}\frac{p_{nm}}{q_n}  \right)&=0\\ \implies z_1+z_2\alpha_1 +z_3\alpha_2+\cdots+z_{m+1}\alpha_m&=0.
\end{align*} 
Hence, $1,\alpha_1,\ldots,\alpha_m$ are $\mathbb{Q}$-linearly independent. This proves the theorem.
\section{Proof of Theorem \ref{maintheorem1}}
We define the  sequences  $(\beta_{1,N})_N$,  $(\beta_{2,N})_N$,  and  $(\beta_{3,N})_N$ of rational numbers as follows.  For each integer $N\geq 1$,  we define 
 $$
 \beta_{1,N}=\sum_{n=1}^{N}\frac{\sigma(n)}{b_1b_2\ldots b_n}=\frac{p_{1,N}}{b_1b_2\ldots b_N}, 
 $$  
 for some positive integer $p_{1,N}$,  
$$
\beta_{2,N}=\sum_{n=1}^{N}\frac{\phi(n)}{b_1b_2\ldots b_n}=\frac{p_{2,N}}{b_1b_2\ldots b_N},
$$ 
for some positive integer $p_{2,N}$,
and
$$
\beta_{3,N}=\sum_{n=1}^{N}\frac{d_n}{b_1b_2\ldots b_n}=\frac{p_{3,N}}{b_1b_2\ldots b_N}.
$$ 
By \eqref{eq1.6}, \eqref{eq1.7} and \eqref{eq1.8}, we easily get
$$
\left|\beta_i-\frac{p_{i,N}}{b_1 b_2 \cdots b_N}\right|<\frac{1}{(b_1 b_2\cdots b_N)^{1+\delta'}};
$$
for infinitely many values of $N\in T$ and for some $\delta'>\frac{1}{3}$. 
\bigskip

By taking $\alpha_i=\beta_i$, $\alpha_2=\alpha'$, $\alpha_3=\alpha''$ and for infinitely many $N\in T$, $p_{iN}=p_{i,N}$ and  $q_N=b_1 b_2 \cdots b_N,$ for $1\leq i \leq 3$  in Theorem \ref{maintheorem}, we get either  $1,\alpha$,  $\alpha'$ and $\alpha''$ are $\mathbb{Q}$-linearly dependent or at least one of them is transcendental. By Theorem \ref{erdos}, we know that $1,\beta_1$, $\beta_2$  and $\beta_3$ are $\mathbb{Q}$-linearly independent. Therefore, we conclude that one of $\beta_1$, $\beta_2$ and  $\beta_3$  is transcendental. This proves the assertation. 
\section{Proof of Theorem \ref{maintheorem2}}
For each integer $1\leq i \leq m$, we define the  sequence  $(\beta_{i,N})_N$  of rational numbers as follows.  For each integer $N\geq 1$,  we define 
 $$
 \beta_{i,N}=\sum_{n=1}^{N}\frac{c_{i,n}}{b_n}=\frac{p_{i,N}}{b_1b_2\ldots b_N},
 $$  
 for some positive integer $p_{i,N}$.
From \eqref{eq1.2},  we obtain that  for each real number $A>1$ there exists a positive constant $N_0$  such that    for all positive integer $N>N_0$ 
$$
\frac{1}{A}.\frac{c_{i,N}}{b_N}>\frac{c_{i,N+1}}{b_{N+1}}.
$$
By this and using mathematical induction we get for every $N$ with $N>N_0$
$$
\frac{1}{A^p}.\frac{c_{i,N}}{b_N}>\frac{c_{i,N+p}}{b_{N+p}}.
$$ 
This implies that for all sufficiently large positive integers $N$, 
\begin{eqnarray*}
\left |\beta_i-\frac{p_{i,N}}{b_1b_2\ldots b_N}\right|&=&\left|\sum_{n=1}^\infty\frac{c_{i,n}}{b_n}-\sum_{n=1}^N\frac{c_{i,n}}{b_n}  \right|=\left|\sum_{n=N+1}^\infty\frac{c_{i,n}}{b_n}\right|\\
&=&\left(\frac{c_{i,N+1}}{b_{N+1}}+\frac{c_{i,N+2}}{b_{N+2}}+\cdots \right)\\
&<&\frac{c_{i,N+1}}{b_{N+1}}\left(1+\frac{1}{A}+\frac{1}{A^2}+\cdots\right)\\
&=&\frac{c_{i,N+1}}{b_{N+1}}\frac{A}{A-1}.
\end{eqnarray*}
Choose $M>\frac{A}{A-1}$. Then by \eqref{eq1.1}, we get infinitely many integers $N$ such that 
$$
\frac{1}{M(b_1 b_2\ldots b_N)^{1+\delta}}>\frac{c_{i,N+1}}{b_{N+1}}.
$$
This implies that  
  $$
\left|\beta_i-\frac{p_{i,N}}{b_1b_2\ldots b_N}\right| < \frac{c_{i,N+1}}{b_{N+1}}\frac{A}{A-1}\leq \frac{1}{(b_1 b_2\ldots b_N)^{1+\delta}}
$$
holds  for infinitely many positive integers $N$. 
\smallskip

By taking $\alpha_i=\beta_i$ and $p_{in}=p_{i,n}$ for $1\leq i \leq m$ in Theorem \ref{maintheorem}, we get that either $1,\beta_1, \beta_2, \ldots,\beta_m$ are $\mathbb{Q}$-linearly dependent or at least one $\beta_i$'s is transcendental. 

\section{proof of Theorem \ref{maintheorem3}}
By Theorem \ref{maintheorem2}, we get that either $1,\beta_1, \beta_2, \ldots,\beta_m$ are $\mathbb{Q}$-linearly dependent or at least one $\beta_i$'s is transcendental. Since the sequences $(c_{i,n})_n$  and $(b_n)_n$ satisfy the hypothesis of Theorem \ref{maintheorem7}, we otain that $1,\beta_1, \beta_2, \ldots,\beta_m$ are $\mathbb{Q}$-linearly independent. Therefore, we conclude that at least one of $\beta_i$'s is transcendental. This proves the theorem. 
\section{Proof of Theorem \ref{maintheorem4}}
For each integer $1\leq i \leq m$, we define the  sequence  $(\beta_{i,N})_N$  of rational numbers as follows.  For each integer $N\geq 1$,  we define 
 $$
 \beta_{i,N}=\sum_{n=1}^{N}\frac{c_{i,n}}{b_n}=\frac{p_{i,N}}{b_1b_2\ldots b_N},
 $$  
 for some positive integer $p_{i,N}$.
From \eqref{eq1.4} we have 
$$
\sqrt[1+\epsilon]{\frac{b_{n+1}}{c_{i,n+1}}}\geq \sqrt[1+\epsilon]{\frac{b_n}{c_{i,n}}}+1.
$$
By this and using mathematical induction we get for all sufficiently large integers $N$ and every integer $r$
$$
\sqrt[1+\epsilon]{\frac{b_{N+r}}{c_{i,N+r}}}\geq \sqrt[1+\epsilon]{\frac{b_N}{c_{i,N}}}+r.
$$
Thus 
\begin{equation*}\label{eq6.1}
\tag{6.1}
\frac{b_{N+r}}{c_{i,N+r}}\geq \left(\sqrt[1+\epsilon]{\frac{b_N}{c_{i,N}}}+r \right)^{1+\epsilon}.
\end{equation*}
Also we have for all real $x>1$
\begin{equation*}\label{eq6.2}
\tag{6.2}
\sum_{s=0}^\infty\frac{1}{(x+s)^{1+\epsilon}}<\int_{x-1}^\infty\frac{dy}{y^{1+\epsilon}}=\frac{1}{\epsilon(x-1)^\epsilon}.
\end{equation*}
From \eqref{eq6.1} and \eqref{eq6.2},  we get that for infinitely many $N$
\begin{eqnarray*}
\left |\beta_i-\frac{p_{i,N}}{b_1b_2\ldots b_N}\right|&=&\left|\sum_{n=1}^\infty\frac{c_{i,n}}{b_n}-\sum_{n=1}^N\frac{c_{i,n}}{b_n}  \right|=\left|\sum_{n=N+1}^\infty\frac{c_{i,n}}{b_n}\right|\\
&=&\left(\frac{c_{i,N+1}}{b_{N+1}}+\frac{c_{i,N+2}}{b_{N+2}}+\cdots \right)\\
&\leq & \left(\sqrt[1+\epsilon]{\frac{b_{N+1}}{c_{i,N+1}}}\right)^{-(1+\epsilon)}+\left(\sqrt[1+\epsilon]{\frac{b_{N+1}}{c_{i,N+1}}}+1\right)^{-(1+\epsilon)}+\cdots\\
&<&\frac{1}{\epsilon}\left(\sqrt[1+\epsilon]{\frac{b_{N+1}}{c_{i,N+1}}}-1\right)^{-\epsilon}.
\end{eqnarray*}
Inequality \eqref{eq1.4} implies that $\lim_{n\rightarrow\infty}\left(\frac{b_n}{c_{i,n}}\right)=\infty$. Hence,  there exists a positive constant $C$ which does not depend on $N$ such that
\begin{eqnarray*}
\left |\beta_i-\frac{p_{i,N}}{b_1b_2\ldots b_N}\right|&<&\frac{1}{\epsilon}\left(\sqrt[1+\epsilon]{\frac{b_{N+1}}{c_{i,N+1}}}-1\right)^{-\epsilon}<\frac{C}{\epsilon}\left(\sqrt[1+\epsilon]{\frac{b_{N+1}}{c_{i,N+1}}}\right)^{\epsilon}\\
&=& \frac{C}{\epsilon}\left(\frac{c_{i,N+1}}{b_{N+1}}\right)^{\frac{\epsilon}{1+\epsilon}}.
\end{eqnarray*}
Choose $M>\frac{C}{\epsilon}$. Then by \eqref{eq1.2}, we get infinitely many integers $N$ such that 
$$
\frac{1}{M(b_1 b_2\ldots b_N)^{1+\delta+\frac{1}{\epsilon}}}>\frac{c_{i,N+1}}{b_{N+1}}.
$$
This implies that  
$$
\left|\beta_i-\frac{p_{i,N}}{b_1b_2\ldots b_N}\right| <\frac{1}{(b_1 b_2\ldots b_n)^{1+\frac{\delta\epsilon}{1+\epsilon}}}
$$
holds  for infinitely many positive integers $N$.    The Rest of the proof is similar to the proof of Theorem \ref{maintheorem3}.
\section{Proof of Theorem \ref{maintheorem5}}
By Theorem \ref{maintheorem4}, we get that either $1,\beta_1, \beta_2, \ldots,\beta_m$ are $\mathbb{Q}$-linearly dependent or at least one $\beta_i$'s is transcendental. Since the sequences $(c_{i,n})_n$  and $(b_n)_n$ satisfy the hypothesis of Theorem \ref{maintheorem7}, we obtain that $1,\beta_1, \beta_2, \ldots,\beta_m$ are $\mathbb{Q}$-linearly independent. Therefore, we conclude that at least one of $\beta_i$'s is transcendental. This proves the theorem. 
\section{Proof of Theorem \ref{maintheorem6}}
For each integer $1\leq i \leq m$, we define the  sequence  $(\beta_{i,N})_N$  of rational numbers as follows.  For each integer $N\geq 1$,  we define 
 $$
 \beta_{i,N}=\prod_{n=1}^{N}\left(1+\frac{c_{i,n}}{b_n}\right)=\frac{p_{i,N}}{b_1b_2\ldots b_N},
 $$  
 for some positive integer $p_{i,N}$.  Consider
\begin{align*}\label{eq9.1}
\left|\beta_i-\frac{p_{i,N}}{b_1 b_2 \ldots b_N}\right|&=\prod_{n=1}^\infty\left(1+\frac{c_{i,n}}{b_n}\right)-\prod_{n=1}^N\left(1+\frac{c_{i,n}}{b_n}\right)\\
&=\prod_{n=1}^N\left(1+\frac{c_{i,n}}{b_n}\right)\left|1-\prod_{n=N+1}^\infty\left(1+\frac{c_{i,n}}{b_n}\right)\right|.\\
\tag{9.1}
\end{align*}
By the hypothesis, for all sufficiently large values of $N$, we have  
$$
\prod_{n=N+1}^\infty\left(1+\frac{c_{i,n}}{b_n}\right)<2\sum_{n=N+1}^\infty \frac{c_{i,n}}{b_n}. 
$$
This implies that 
$$
\left|1-\prod_{n=N+1}^\infty\left(1+\frac{c_{i,n}}{b_n}\right)\right|< 3\sum_{n=N+1}^\infty\frac{c_{i,n}}{b_n}.
$$
Thus, by \eqref{eq9.1} we have
\begin{equation*}
\left|\beta_i-\frac{p_{i,N}}{b_1 b_2 \ldots b_N}\right|\leq 3\prod_{n=1}^N\left(1+\frac{c_{i,n}}{b_n}\right)\left(\sum_{n=N+1}^\infty\frac{c_{i,n}}{b_n}\right).
\end{equation*}
Using a similar argument such as Theorem \ref{maintheorem3}, from \eqref{eq1.2},  we conclude that for all sufficiently large positive integers $N$, 
\begin{eqnarray*}
\left|\beta_i-\frac{p_{i,N}}{b_1 b_2 \ldots b_N}\right|\leq 3\prod_{n=1}^N\left(1+\frac{c_{i,n}}{b_n}\right)\left(\sum_{n=N+1}^\infty\frac{c_{i,n}}{b_n}\right)<3\prod_{n=1}^N\left(1+\frac{c_{i,n}}{b_n}\right)\frac{c_{i,N+1}}{b_{N+1}}\frac{A}{A-1}.
\end{eqnarray*}  
Hence, by \eqref{eq1.1} we obtain
\begin{equation*}\label{eq9.2}
\tag{9.2}
\left|\beta_i-\frac{p_{i,N}}{b_1 b_2 \ldots b_N}\right|< \prod_{n=1}^N\left(1+\frac{c_{i,n}}{b_n}\right)\frac{1}{(b_1b_2\ldots b_N)^{1+\delta}},
\end{equation*}
for infinitely many values of $N$. 
\smallskip

By the assumption $\frac{c_{i,n}}{b_n}\leq 1$, for $n\geq 1$, we have 
\begin{equation*}
\prod_{n=1}^N\left(1+\frac{c_{i,n}}{b_n}\right)<2^N
\end{equation*}
for all integer $N\geq 1$.  Thus from \eqref{eq9.2}, we have
$$
\left|\beta_i-\frac{p_{i,N}}{b_1 b_2 \ldots b_N}\right|<\frac{2^N}{(b_1 b_2\ldots b_N)^{1+\delta}}.
$$
Since the sequence $(b_n)_n$ grows like a doubly exponential sequence, we can find $\frac{1}{m}<\delta'<\delta$ such that 
$$
\frac{2^N}{(b_1 b_2 \ldots b_N)^{1+\delta}}<\frac{1}{(b_1 b_2 \ldots b_N)^{1+\delta'}}.
$$
This implies that for each $1\leq i \leq m$
$$
\left|\beta_i-\frac{p_{i,N}}{b_1 b_2 \ldots b_N}\right|<\frac{1}{(b_1 b_2 \ldots b_N)^{1+\delta'}}
$$
 holds for infinitely many valus of $N$.  By taking $\alpha_i=\beta_i$ and $p_{in}=p_{i,n}$ for $1\leq i \leq m$ in Theorem \ref{maintheorem}, we get that either $1,\beta_1, \beta_2, \ldots,\beta_m$ are $\mathbb{Q}$-linearly dependent or at least one $\beta_i$'s is transcendental. Since the sequences $(c_{i,n})_n$  and $(b_n)_n$ satisfy the hypothesis of Theorem \ref{maintheorem8}, we obtain that $1,\beta_1, \beta_2, \ldots,\beta_m$ are $\mathbb{Q}$-linearly independent. Therefore, we conclude that at least one of $\beta_i$'s is transcendental. This proves the theorem. 
\section{Proof of Corollary \ref{cor}}
By taking $c_{1,n}=1$ and $c_{2,n}=d(n)$, we see that these sequences satisfy the hypothesis of Theorem \ref{maintheorem2}. Hence, by Theorem \ref{maintheorem2}, we get that  either $1, \alpha$ and $\alpha'$ are $\mathbb{Q}$-linearly dependent or at least one of $\alpha$ and $\alpha'$ is transcendental.  In order to finish the proof of this corollary, we shall prove that $1, \alpha$ and $\alpha'$ are $\mathbb{Q}$-linearly independent. 
\smallskip

Suppose that these numbers are $\mathbb{Q}$-linearly dependent. Then,  there exist  integers $z_0, z_1$ and $z_2$ not all zero such that 
\begin{equation*}\label{eq10.1}
\tag{10.1}
z_0+ z_1 \sum_{n=1}^\infty\frac{1}{b_n}+ z_2 \sum_{n=1}^{\infty}\frac{d(n)}{b_n}=0.
\end{equation*}
This is equivalent to 
$$
z_0 + z_1 \sum_{n=1}^N\frac{1}{b_n} + z_2  \sum_{n=1}^N\frac{d(n)}{b_n}=-\left(z_1 \sum_{n=N+1}^\infty\frac{1}{b_n}+ z_2 \sum_{n=N+1}^\infty\frac{d(n)}{b_n} \right).
$$
By multiplying $b_1 b_2 \cdots b_N$ on both sides, we get
$$
b_1 b_2\cdots b_N\left( z_0 + z_1 \sum_{n=1}^N\frac{1}{b_n} + z_2  \sum_{n=1}^N\frac{d(n)}{b_n}\right)=-b_1  \cdots b_N\left(\sum_{n=N+1}^\infty\frac{z_1}{b_n}+  \sum_{n=N+1}^\infty\frac{ z_2 d(n)}{b_n} \right).
$$
We note that the left-hand side of this is an integer.  Now,  we claim the following.
\smallskip

\noindent{\bf Claim.} The quantity 
$$
\left|-b_1 b_2 \cdots b_N\left(z_1 \sum_{n=N+1}^\infty\frac{1}{b_n}+ z_2 \sum_{n=N+1}^\infty\frac{d(n)}{b_n} \right)\right|\rightarrow 0 \quad \mbox{as} \quad N\rightarrow\infty.
$$
In order to prove the claim, we estimate the above quantity as follows.   Consider
\begin{align*}
\left|-b_1 b_2 \cdots b_N\left(z_2 \sum_{n=N+1}^\infty\frac{d(n)}{b_n}\right)\right| &\leq |z_2|\left(\frac{d(N+1)}{b_{N+1}}+\frac{d(N+2)}{b_{N+2}}+\cdots\right)\\
\end{align*}
Using $d(n)=O(n)$, we have
\begin{align*}
\left|-b_1 b_2 \cdots b_N\left(z_2 \sum_{n=N+1}^\infty\frac{d(n)}{b_n}\right)\right| &\leq |z_2|\left(\frac{d(N+1)}{b_{N+1}}+\frac{d(N+2)}{b_{N+2}}+\cdots\right)\\
&<\frac{1}{b_1 b_2 \cdots b_N}\left(\frac{N+1}{b_1 b_2 \cdots b_N}+\frac{N+2}{(b_1 b_2 \cdots b_N)^{6N}+\cdots}\right)\\
&<\frac{C}{b_1 b_2 \cdots b_N}.
\end{align*}
Hence,
\begin{equation*}\label{eq10.2}
\tag{10.2}
\left|-b_1 b_2 \cdots b_N\left(z_2 \sum_{n=N+1}^\infty\frac{d(n)}{b_n}\right)\right|\rightarrow 0  \quad \mbox{as} \quad N\rightarrow\infty.
\end{equation*}
Similarly we get
\begin{equation*}\label{eq10.3}
\tag{10.3} 
\left|-b_1 b_2 \cdots b_N\left(z_1 \sum_{n=N+1}^\infty\frac{1}{b_n}\right)\right|\rightarrow 0 \quad \mbox{as} \quad N\rightarrow\infty.
\end{equation*}
Thus, by \eqref{eq10.2} and \eqref{eq10.3}, we get the claim.
\bigskip
Hence, we have  
\begin{equation*}\label{eq10.4}
\tag{10.4}
 z_0 + z_1 \sum_{n=1}^N\frac{1}{b_n} + z_2  \sum_{n=1}^N\frac{d(n)}{b_n}=0 
\end{equation*}
for all sufficiently large values of $N$.  By multiplying $b_1 b_2 \cdots b_{N-1}$, we get
\begin{align*}
b_1 b_2\cdots b_{N-1}\left(z_0 + z_1 \sum_{n=1}^{N-1}\frac{1}{b_{n-1}} + z_2  \sum_{n=1}^{N-1}\frac{d(n)}{b_{n-1}}\right)&=\frac{-b_1 b_2 \cdots b_{N-1}(z_1 +z_2 d(N)}{b_N}\\
&=\frac{z_1 +z_2 d(N))}{b_1 b_2 \cdots b_{N-1}}.
\end{align*}
Clearly, the left hand side is an integer.  
Since $\frac{(|z_1|+|z_2|d(N))}{b_1 b_2 \cdots b_{N-1}}\rightarrow 0$ as $N\rightarrow \infty$, we get
$$
0\leq \left|\frac{(|z_1|+|z_2|d(N))}{b_1 b_2 \cdots b_{N-1}}\right|<1
$$
for all sufficiently large values of $N$. Thus, we have 
\begin{equation*}\label{eq10.5}
\tag{10.5}
z_0 + z_1 \sum_{n=1}^{N-1}\frac{1}{b_n} + z_2  \sum_{n=1}^{N-1}\frac{d(n)}{b_n}=0.
\end{equation*}
Hence, by \eqref{eq10.4} and \eqref{eq10.5},  we get 
$$
\frac{z_1  + z_2 d(N)}{b_N}=0\iff \frac{1}{d(N)}=-\frac{z_2}{z_1}
$$
for all sufficiently large values of $N$.   This implies that the sequence $(d(n))_n$ is eventually constant, which is a contradiction of the fact that it has at least two limit points.

\end{document}